\documentclass[11pt, a4paper]{amsart}

\usepackage{amsmath,amssymb,amsthm}
\numberwithin{equation}{section}
\newtheorem{theorem}{Theorem}[section]

\newtheorem{lemma}[theorem]{Lemma}
\theoremstyle{remark}
\newtheorem{remark}{Remark}[section]

\theoremstyle{definition}
\newtheorem{definition}[theorem]{Definition}

\newcommand{\R}{\mathbb{R}}
\newcommand{\C}{\mathbb{C}}

\begin{document}

%%%%TITLE%%%%%

\title[Helmholtz equation in exterior domains]
{A priori estimates for the Helmholtz equation with electromagnetic potentials in exterior domains}

%%%%ABSTRACT%%%%

\begin{abstract}
  We study the Helmholtz equation with electromagnetic-type perturbations, in the exterior of a domain, in dimension $n\geq3$. We prove, by multiplier techniques in the sense of Morawetz, a family of a priori estimates from which the limiting absorption principle follows. Moreover, we give some standard applications to the absence of embedded eigenvalues and zero-resonances, under explicit conditions on the potentials.
\end{abstract}

\date{\today}

%%%%AUTHORS%%%%

\author{Juan Antonio Barcel\'o}
\address{Juan Antonio Barcel\'o: ETSI de Caminos, Universidad Polit\'ecnica de Madrid, 28040, Madrid,
Spain}
\email{ juanantonio.barcelo@upm.es}

\author{Luca Fanelli}
\address{Luca Fanelli: Universidad del Pais Vasco, Departamento de
Matem\'aticas, Apartado 644, 48080, Bilbao, Spain}
\email{luca.fanelli@ehu.es}

\author{Alberto Ruiz}
\address{Alberto Ruiz: Departamento de Matem\'aticas, Universidad Aut\'onoma de Madrid, 28049, Madrid,
Spain}
\email{alberto.ruiz@uam.es}

\author{Maricruz Vilela}
\address{Maricruz Vilela: Departamento de Matem\'atica Aplicada, Universidad de Valladolid, Plaza
Santa Eulalia 9 y 11, Segovia, Spain}
\email{maricruz@dali.eis.uva.es}

%%%%SUBJCLASS & KEYWORDS%%%%



\maketitle

%%%% BODY OF PAPER %%%%

\section{Introduction}\label{sec:intro}

We study of the following Helmholtz equation
\begin{equation}\label{eq:helmholtz-laplace}
  (\Delta_A-V(x))u+(k^2\pm i\epsilon)u=f
  \quad
  \text{in }E,
\end{equation}
where $k\in\R$, $E=\R^n\setminus\Omega$ and $\Omega\subset\R^n$ is a bounded domain; moreover, we denote by
\begin{equation*}
  \nabla_A=\nabla-iA,
  \qquad
  \Delta_A=\nabla_A^2,
\end{equation*}
with $A:\R^n\to\R^n$, $V:\R^n\to\R$. When $A\equiv0$, \eqref{eq:helmholtz-laplace} is the
 standard Helmholtz equation, where the coefficient $V$ represents an external potential. The operator $H=-\Delta_A+V$ standardly represents an electromagnetic Schr\"odinger hamiltonian,
   which in quantum mechanical models describes the interactions of free particles with some external electromagnetic fields.
    In this paper, we are interested in the so called Agmon-H\"ormander estimates, which cover the full frequency range $k^2>0$. The seminal papers by Agmon and H\"ormander \cite{a}, \cite{a2}, \cite{ah}, in the free case $A\equiv V\equiv0$, inspired a huge literature, which has been produced in order to obtain weighted $L^2$ estimates for solutions of Helmholtz equations. 
As it is well known, one of the consequences of the Agmon-H\"ormander estimate is related to the singular spectrum of the operator $H$.    
    In what we call the
{\it electric} case $A\equiv 0$,
Agmon \cite{a} and Saito \cite{s} proved that for
potentials $V(x)$ which decay at infinity as $V(x)=O(|x|^{-1})$,
the singular spectrum is absent in $(0, \infty)$ for
the operator $H=-\Delta +V$ (see also \cite{brv0}). Later, Lavine \cite{l} and Arai \cite{ar} studied the same problem for {\it repulsive} potentials, i.e. $\partial V / \partial |x| \leq
0$. 

In all these cases, Fourier Analysis is involved as a
crucial tool in the proofs strategy; on the other hand, the Fourier transform does not permit in general to treat neither rough potentials neither the case in which the same problems are settled in domains which are different from the whole space. For this reason, a great effort has been spent  in order to develop multiplier methods which work directly on the equation, inspired to the techniques introduced by Morawetz in \cite{mor} for the Klein-Gordon equation. For the literature about the Helmholtz equation, we first mention Ikebe and Saito \cite{IS}. In this paper, the authors prove existence and uniqueness of solutions of \eqref{eq:helmholtz-laplace} in the whole space, satisfying some weighted $L^2$-estimates and a suitable Sommerfeld radiation condition. In fact, they consider a compact set of frequencies $k^2\in[k_0,k_1]$, $k_0>0$, which are far away from the origin. This permits them to treat also potentials which decay slowly at infinity.
In the same spirit, using a refinement of the standard Morawetz multipliers, Perthame and Vega in \cite{pv2} detected the Sommerfeld radiation condition, in the case $A\equiv0$, and for the complete set of large frequencies $k^2>k_0>0$, for potentials which can vary at infinity. 
In a previous work \cite{pv}, the same authors had developed the Morawetz technique in orther to recover the Agmon-H\"ormander estimate, in Morrey-Campanato spaces, in the full frequency-range $k^2>0$. We also mention \cite{brvv}, in which the authors relate the problem with the weak dispersive properties of the Schr\"odinger evolution equation, and can also treat long range potentials $V$.
In the same spirit as in \cite{pv}, Fanelli in \cite{F} extended the argument by Perthame and Vega to the {\it magnetic} case $A\neq0$ (see also \cite{F2} and the references therein, for a survey about the topic).

As it is well known, the fact that the Agmon-H\"ormander estimate holds for all $k^2>0$ requires that the potentials involved are short-range (see the rest of the paper); in fact this estimate has several consequences about weak dispersive properties of the Schr\"odinger flow $e^{itH}$, as local smoothing estimates, which in general do not hold with long-range potentials (see among the others \cite{brv}, \cite{pda-lf1}, \cite{FV}, \cite{k}, \cite{ky}, \cite{lp}, \cite{ls}).
The aim of this paper is to prove, in the same spirit as in \cite{pv} and \cite{F}, the same a priori estimates, for all positive frequencies $k^2>0$, for solutions of equation \eqref{eq:helmholtz-laplace} in the exterior of a domain, with Dirichlet or Neumann boundary data.

Before stating our main results, we need some preliminary definitions. The magnetic field $B$ is the anti-symmetric gradient of the field $A$, namely
\begin{equation*}
  B=DA-(DA)^t,
  \qquad
  (DA)_{ij}=\frac{\partial A^i}{\partial x_j},
  \quad
  (DA)^t_{ij}=(DA)_{ji}.
\end{equation*}
In geometrical terms, $B$ is nothing else than the differential of the linear 1-form $\omega=A^1(x)dx_1+
\dots+A^n_xd_{x_n}$ which is naturally associated to $A$, i.e. $B=d\omega$. In particular, in dimension $n=3$, the magnetic field $B$ is identified as
$B=\text{curl}A$, due to the isomorphism between 1-forms and
2-forms; this fact has to be interpreted in terms of the action
\begin{equation*}
  Bv=\text{curl}A\times v,
  \qquad
  \text{for all }v\in\R^3,
\end{equation*}
where the cross is the vectorial product on $\R^3$.
Following \cite{FV}, we denote by
$B_\tau:\R^n\to\R^n$ the tangential component of the magnetic
  field $B$, given by
  \begin{equation}\label{def.bitau}
    B_\tau(x):=\frac{x}{|x|}B.
  \end{equation}
  Hence the $i$-th component $B_\tau^i$ of the vector $B_\tau$ is given by
  \begin{equation*}
    B_\tau^i=\sum_{j=1}^n\frac{x_j}{|x|}\left(A^i_j-A^j_i\right),
    \qquad
    A^i_j:=\frac{\partial A^i}{\partial x_j}.
  \end{equation*}
  Observe that in dimension $n=3$ it coincides with
  \begin{equation*}
    B_\tau(x):=\frac{x}{|x|}\times\text{curl}A(x),
  \end{equation*}
  the cross denoting the vectorial product in $\R^3$.
  In addition, we say that $B$ is \textit{non-trapping} if $B_\tau=0$.

As a $3D$-example
of potential $A$ for which $B_\tau\equiv0$, take
  \begin{equation}\label{eq.example}
    A=\frac1{x^2+y^2+z^2}(-y,x,0)=\frac{1}{x^2+y^2+z^2}(x,y,z)\wedge(0,0,1).
  \end{equation}
  We can check that
  \begin{equation*}
    \nabla\cdot A=0,
    \qquad
    B=-2\frac{z}{(x^2+y^2+z^2)^2}(x,y,z),
    \qquad
    B_\tau=0.
  \end{equation*}
Another (more singular) example is the following:
  \begin{equation}\label{eq.example2}
    A=\left(\frac{-y}{x^2+y^2},\frac{x}{x^2+y^2},0\right)=
    \frac{1}{x^2+y^2}(x,y,z)\wedge(0,0,1).
  \end{equation}
  Here we have $B=(0,0,\delta)$, with $\delta$ denoting Dirac's delta function. Again we have $B_\tau=0$ . By translations, we can produce the same kind of examples with a singularity at a generic point $x_0\in R^n$. For a larger class of examples, see \cite{FV}.
  
We now pass to introduce the abstract functional setting in which we work in the sequel. In order to do this, we first need to introduce some regularity assumptions on the Hamiltonian $H$.

{\bf (H1)} The Hamiltonian $H=-\Delta_A+V$ is self-adjoint (and positive) on $L^2(\R^n)$, with form domain
\begin{equation*}
  \mathcal D(H)=\{f\in L^2:\int_{\R^n}\left(|\nabla_Af|^2+V|f|^2\right)<\infty\}.
\end{equation*}  
Assumption (H1) has several consequences about the existence theory for equation 
\eqref{eq:helmholtz-laplace}. For the kind of potentials we deal with in the sequel, this can be standardly proved by perturbation theory, under suitable conditions on the potentials. Indeed, one could argue by first proving that $\Delta_A$ is self-adjoint, under local integrability conditions on $A$, and then assume that the negative part of $V$ is a perturbation of $\Delta_A$ in the Kato-Rellich sense. We prefer here to state (H1) as an abstract requirement (see the standard references \cite{CFKS}, \cite{LS} for details).

One of the consequences of assumption (H1) is that, via Spectral Theorem, we can define the positive powers $H^s$ of the operator $H$, and the distorted Sobolev norms
\begin{equation}\label{eq:sobolev}
  \|f\|_{\mathcal H^s}:=\|H^{\frac s2}f\|_{L^2},
  \qquad
  s\geq0;
\end{equation}
these are the natural spaces in which we work below.
We also need to assume the following:

{\bf (H2)}  the Sobolev spaces $H^s$ and $\mathcal H^s$ are equivalent.

Assumption (H2) will be used in the sequel to justify a Trace Theorem for the spaces $\mathcal H^s$, by means of the usual one for the case $A,V=0$. Conditions on $A,V$ for the validity of (H2) can be found in Theorem 1.2 in \cite{DFVV}, case $q=2$. In what follows, (H1)-(H2) will be always implicitly assumed.

We finally introduce some notations. In the following, we always
denote by $B_y(R)=\{x\in\R^n:|x-y|\leq R\}$ and $S_y(R)=\partial B_y(R)$; we also use the notations $B(R)=B_0(R)$, $S(R)=S_0(R)$, $C(j)=C_0(j)$, $N(f)=N_0(f)$.

Given a set $D\subset\R$ and a function $f:D\to\R$, we write $f=f_+-f_-$ as the difference of the positive and the negative part, respectively, of $f$.
Moreover, let $D\subset\R^n$, $f,g:D\to\C$ be Borel-measurable functions and identify $f,g$ with their trivial extension to $\R^n$; for any $y\in\R^n$ and $p\geq1$, we define
\begin{equation}\label{eq:spazi}
  \|f\|_{L^p_rL^\infty(S_y(r))}=
  \left(\int_0^{\infty}
  \sup_{|x-y|=\rho}|f(x)|^pd\rho\right)^{\frac1p},
\end{equation}
Analogously, let $C_y(j)=\{x\in\R^n:2^j\leq|x-y|<2^{j+1}\}$,
and denote by
\begin{equation}\label{eq:N}
  N_y(f)=\sum_{j\in\mathbb Z}\left(2^{j+1}\int_{C_y(j)}|f|^2
  \right)^{\frac12};
\end{equation}
we easily notice the duality relation
\begin{equation*}
  \int_D fg\,dx
  \leq
  N_y(f)\cdot\left(\sup_{R>0}\frac1R
  \int_{D\cap(B_y(R))}|g|^2dx
  \right)^{\frac12}.
\end{equation*}

Our first result concerns with equation \eqref{eq:helmholtz-laplace} with Dirichlet boundary conditions in an exterior domain $E$, in dimension $n\geq3$.
\begin{theorem}\label{thm:lap3}
  Let $n\geq3$; let $\Omega\subset\R^n$ be a bounded domain with Lipshitz boundary and $E=\R^n\setminus \Omega$ be its exterior.
  Let $y\in\R^n\setminus\partial E$, and denote by
  \begin{equation}\label{eq:convexity}
    \partial E^-_{y}=\{x\in\partial E:\frac{x-y}{|x-y|}\cdot\eta\leq 0\},
    \qquad
    \partial E^+_{y}=\{x\in\partial E:\frac{x-y}{|x-y|}\cdot\eta\geq 0\},
  \end{equation}
  where $\eta(x)$ denotes the unit outgoing normal to the boundary $\partial E$. Let
  \begin{equation}\label{eq:hardy}
    V=V_+-V_-,
    \qquad
    V_-\in L^\infty(E);
  \end{equation}
  assume moreover that
  \begin{align}\label{eq:smallness3}
    &
    \||\cdot-y|^{\frac32}B_\tau\|_{L^2_rL^\infty(S_y(r))}
    +\||\cdot-y|^2(\partial_rV)_+\|_{L^1_rL^\infty(S_y(r))}
    \\
    &
    +\|\,|\cdot-y|V_+\|_{L^1_rL^\infty(S_y(r))}
    <\delta,
    \nonumber
  \end{align}
  in dimension $n=3$, or
  \begin{equation}\label{eq:smallness4}
    \||\cdot-y|^2B_\tau\|_{L^\infty}
    +\||\cdot-y|^3(\partial_rV)_+\|_{L^\infty}
    +\|\,|\cdot-y|^2V_+\|_{L^\infty}<\delta,
  \end{equation}
  in dimension $n\geq4$, for a sufficiently small $\delta>0$.
  Then any solution $u\in\mathcal H^1(E)\cap\mathcal H^{\frac32+}_{loc}$, close to $\partial E$, of
  \begin{equation}\label{eq:diric-laplace}
    \begin{cases}
      (\Delta_A-V(x))u+(k^2\pm i\epsilon)u=f
      \quad
      \text{in }E
      \\
      u\vert_{\partial E}\equiv0
    \end{cases}
  \end{equation}
  satisfies
  \begin{align}\label{eq:estimate1}
    &
    \int_E\frac{|\nabla_A^\bot u|^2}{|x-y|}+
    \sup_{R>0}\left(
    \frac1R\int_{E\cap B_y(R)}\left(|\nabla_Au|^2
    +k^2|u|^2\right)
    +\frac1{R^2}\int_{E\cap S_y(R)}|u|^2\right)
    \\
    & \ \ \ 
    +(n-3)\int_E\frac{|u|^2}{|x-y|^3}
    +\sup_{R>0}\frac1R\int_{E\cap B_y(R)}V_-|u|^2
    +\int_{E}(\partial_rV)_-|u|^2
    \nonumber
    \\
    & \ \ \ 
    +\int_{\partial E^-_y}|\nabla_A^\eta u|^2\left(-\frac{x-y}{|x-y|}\cdot\eta\right)\,d\sigma
    \nonumber
    \\
    &
    \leq
    C\left[N_y(f)^2+(|\epsilon|+k^2+\|V_-\|
    _\infty)\left(
    N_y\left(\frac{f}{|k|}\right)\right)^2\right]
    \nonumber
    \\
    & \ \ \
    +C\int_{\partial E^+_y}|\nabla_A^\eta u|^2
    \left(\frac{x-y}{|x-y|}\cdot\eta\right)\,d\sigma,
    \nonumber
  \end{align}
  for some $C>0$, where $\nabla_A^\eta u=\nabla_Au\cdot\eta$, and $|\nabla_A^\bot u|^2 = |\nabla_Au|^2-\left|\nabla_Au\cdot\frac{x-y}{|x-y|}\right|^2$
  is the tangential component of $\nabla_Au$ to the sphere $S_y(1)$.
  \end{theorem}
  \begin{remark}
    The regularity assumption $u\in\mathcal H^1(E)$ is needed to justify the terms in estimate \eqref{eq:estimate1} involving integrals in $E$. Notice that, for $f\in L^2$, it follows by standard elliptic theory. On the other hand, the requirement $u\in\mathcal H^{\frac32+}_{loc}$, together with assumption (H2), justifies all the boundary terms in the estimate, via the usual Trace Theorem.
    \end{remark}
 \begin{remark}
  Notice that, in the case $\Omega=\R^n$, estimate \eqref{eq:estimate1} recovers
   the ones proved in \cite{F}.
\end{remark}
\begin{remark}[3D-estimate]
  We remark that the term containing $\int|u|^2/|x-y|^3$ at the left-hand side of \eqref{eq:estimate1} is not present in the estimates in the case $n=3$ (indeed, the integral term is too singular at $y$); this justify the choice to write the constant $(n-3)$.
\end{remark}
\begin{remark}[Decay at infinity]
  The following potentials  
  \begin{equation*}
    B_\tau,V_+=\frac C{|x|^{2-\epsilon}+|x|^{2+\epsilon}},
    \qquad
    (\partial_rV)_+=\frac C{|x|^{3-\epsilon}+|x|^{3+\epsilon}}
    \qquad
    (n=3)
  \end{equation*}
  \begin{equation*}
    B_\tau,V_+=\frac C{|x|^2},
    \qquad
    (\partial_rV)_+=\frac C{|x|^3}
    \qquad
    (n\geq4),
  \end{equation*}
  for a sufficiently small constant $C>0$, and a small $\epsilon>0$, satisfy assumptions \eqref{eq:smallness3} and \eqref{eq:smallness4} (in the case $y=0$). We stress again that, if one looks to the same estimates for large frequencies $k^2\geq k_0>0$, weaker decay conditions are sufficient (see \cite{IS}, \cite{pv2}).
\end{remark}
\begin{remark}[Limiting absorption principle]
  Theorem \ref{thm:lap3} implies in a standard way (see \cite{a}) the so called {\it limiting absorption principle}, for which the resolvent operator of $H$ can be extended, on the positive real line, to a bounded operator between weighted $L^2$-spaces. Moreover, starting by the same estimate \eqref{eq:estimate1}, it is possible to obtain the appropriate Sommerfeld radiation condition implying uniqueness of solutions, as it was recently proved by Zubeldia in \cite{Z}.
\end{remark}
\begin{remark}[Absence of zero-resonances and local smoothing]
  Consider estimate \eqref{eq:estimate1}, in the case of Dirichlet boundary conditions. In particular, when $\Omega$ is star-shaped with respect to 0, no boundary terms appear at the right-hand side. In particular, it implies the absence of embedded eigenvalues for $H=-\Delta_A+V$ in the positive line $[0,+\infty)$. Moreover, following \cite{brv} and \cite{F}, we give a natural definition of zero-resonance:
  \begin{definition}
    A function $u$ is a zero-resonance for Dirichlet boundary conditions if
    \begin{align*}
      & u\vert_{\partial E\equiv0,
      \qquad
      }u\notin L^2,
      \qquad
      H^{\frac12}u\in L^2_{\text{loc}},
      \qquad
      |V|^{\frac12}u\in L^2
      \\
      &
      \sup_{R>1}\int_{E\cap B_0(R)}|u|^2\left(|V|^2+\frac1{1+|x|^2}\right)<\infty
      \\
      &
      \liminf_{R\to+\infty}\int_{E\cap B_0(R)}|u|^2\left(|V|^2+\frac1{1+|x|^2}\right)=0
    \end{align*}
  \end{definition}  
   One could easily see that estimate \eqref{eq:estimate1} implies that the operator $H$ has no zero-resonances (see also \cite{brv}, \cite{F}). 
 \end{remark}
 \begin{remark}[Local smoothing]
   The local smoothing  property for solutions of dispersive equations, like the Schr\"odinger and the wave equations, is
 a standard application of the Agmon-Hormander estimate for the Helmholtz equation. This was firstly remarked by Kato in \cite{k}. On the other hand, one of the advantages of the multiplier techniques is that they permit to work directly on the equations, without needing an abstract theorem. In the last years, multiplier techniques for the local smoothing have been developed in \cite{brv}, \cite{pda-lf1}, \cite{FV} (see also \cite{F2} for a survey about the topic). The same arguments in all these papers can be repeated in the case of exterior domains, arguing exactly as in the proofs of our theorems. We omit here straightforward details.
\end{remark}
Notice that, if $\Omega$ is star-shaped with respect to $y$, then $\partial E\equiv\partial E^-_y$,
  hence no boundary terms appear at the right-hand side of estimates \eqref{eq:estimate1}.
  In this case, we can in fact prove a more general result, via the same techniques, with generic boundary conditions.
  
  \begin{theorem}\label{thm:main}
  Let $n\geq3$; let $\Omega\subset\R^n$ be a bounded domain with Lipshitz boundary, $E=\R^n\setminus \Omega$ be its exterior, with $0\notin\partial\Omega$. Assume that there exists $\beta>0$ such that 
  $\frac x{|x|}\cdot\eta(x)\leq-\beta$, for any $x\in\partial E$,
  where $\eta(x)$ denotes the unit outgoing normal to the boundary $\partial E$ in $x$. Assume \eqref{eq:hardy}, and moreover that $V,A\in L^\infty(\partial E)$. 
  In addition, assume \eqref{eq:smallness3}, in dimension $n=3$ or \eqref{eq:smallness4} in dimension $n\geq4$ (with $y=0$).
  
  Then any solution $u\in\mathcal H^1(E)\cap\mathcal H^{\frac32+}_{loc}$, close to $\partial E$, of
  \eqref{eq:helmholtz-laplace}  satisfies
  \begin{align}\label{eq:estimate1main}
    &
    \int_E\frac{|\nabla_A^\bot u|^2}{|x|}+
    \sup_{R>0}\left(
    \frac1R\int_{E\cap B(R)}\left(|\nabla_Au|^2
    +k^2|u|^2\right)
    +\frac1{R^2}\int_{E\cap S(R)}|u|^2\right)
    \\
    &
    +(n-3)\int_E\frac{|u|^2}{|x|^3}
    +\sup_{R>0}\frac1R\int_{E\cap B(R)}V_-|u|^2
    +\int_{E}(\partial_rV)_-|u|^2
    \nonumber
    \\
    &
    +\beta\int_{\partial E}|\nabla_A^\eta u|^2\,d\sigma
    \nonumber
    \\
    &
    \ \ \
    \leq
    C\left[N_y(f)^2+(|\epsilon|+k^2+\|V_-\|
    _\infty)\left(
    N_y\left(\frac{f}{|k|}\right)\right)^2\right]
    \nonumber
    \\
    & \ \ \
    +C\left(\frac1\beta+1\right)\int_{\partial E}|u|^2\,d\sigma
    +C\left(\frac1\beta+1\right)\int_{\partial E}|\nabla_A^\tau u|^2\,d\sigma
    \nonumber
  \end{align}
  for some $C>0$, where $\nabla_A^\eta u=\nabla_Au\cdot\eta$, $|\nabla_A^\tau u|^2=|\nabla_Au|^2-|\nabla_A^\eta u|^2$, with $\nabla_A^\eta u\cdot\nabla_A^\tau u=0$ is the tangential component of $\nabla_Au$ to the boundary $\partial E$, and $|\nabla_A^\bot u|^2 = |\nabla_Au|^2-\left|\nabla_Au\cdot\frac{x-y}{|x-y|}\right|^2$
  is the tangential component of $\nabla_Au$ to the sphere $S(1)$.
  
  Moreover, there exists $k_0=k_0(\beta,A,V)>0$ such that, for any $k\geq k_0$, we have
  \begin{align}\label{eq:estimate1main2}
    &
    \int_E\frac{|\nabla_A^\bot u|^2}{|x|}+
    \sup_{R>0}\left(
    \frac1R\int_{E\cap B(R)}\left(|\nabla_Au|^2
    +k^2|u|^2\right)
    +\frac1{R^2}\int_{E\cap S(R)}|u|^2\right)
    \\
    &
    +(n-3)\int_E\frac{|u|^2}{|x|^3}
    +\sup_{R>0}\frac1R\int_{E\cap B(R)}V_-|u|^2
    +\int_{E}(\partial_rV)_-|u|^2
    \nonumber
    \\
    &
    +\beta\int_{\partial E}|\nabla^\eta u|^2\,d\sigma
    +\beta k^2\int_{\partial E}|u|^2\,d\sigma
    \nonumber
    \\
    &
    \ \ \
    \leq
    C\left[N_y(f)^2+(|\epsilon|+k^2+\|V_-\|
    _\infty)\left(
    N_y\left(\frac{f}{|k|}\right)\right)^2\right]
    \nonumber
    \\
    & \ \ \
    +C\left(\frac1\beta+1\right)\int_{\partial E}|\nabla^\tau u|^2\,d\sigma.
    \nonumber
  \end{align}
\end{theorem}
\begin{remark}
  As \eqref{eq:estimate1main2} shows, for high frequencies $k\geq k_0>0$ we have a better estimate; indeed we can also control the boundary term depending on $u$ at the left-hand side of \eqref{eq:estimate1main2}. 
  
  There are several estimates for the radiating solutions of the exterior Dirichlet problem in the literature in the case $A=0$ and $V=0.$
  In particular, Alber (see \cite{alber}) studied the 2D-case (see also \cite {ramm} for $n=3$) and proved that for $\alpha>0,$
 \begin{align*}
 &
  \frac{1}{R^{1+\alpha}}\int_{E\cap B(R)}|\nabla u|^2
  +
   \frac{k^2}{R^{1+\alpha}}\int_{E\cap B(R)}|u|^2
   \\
   &
   \leq
   \int_{E\cap B(R)}|x|^2|f|^2
   +
   \int_{\partial E}|\nabla^{\tau}u|^2
   +
   (1+k^2)\int_{E\cap B(R)}|u|^2.
\end{align*}
Morawetz and Ludwig proved (see \cite{moralu}) in dimension $n=3$ that
\begin{align*}
 &
\beta\int_{\partial E}|\nabla^{\eta} u|^2
  +
\int_{E}|\nabla^{\bot} u|^2
+
\frac{1}{2}\int_{E}\left|\partial_r u-iku+\frac{(n-1)}{2|x|}u\right|^2 
   \\
   &
   \leq
   \int_{E}|x|^2|f|^2
   +
  \frac{C}{\beta} \int_{\partial E}|\nabla^{\tau}u|^2
   +
(1+k^2) \frac{C}{\beta}\int_{\partial E}|u|^2.
\end{align*}
In the 3D-case, estimates \eqref{eq:estimate1main}, \eqref{eq:estimate1main2} are stronger compared with the previous ones, which have to be understood by limiting absorption principles as limits when $\epsilon$ tends to $0$ of solutions of \eqref{eq:helmholtz-laplace}.
Notice that we improve the powers of $R$ and $|x|$ and the dependence on $k.$
In addition, we can also treat electromagnetic perturbations of the Helmholtz equation.
  \end{remark}
The rest of the paper is devoted to the proofs of Theorems \ref{thm:lap3} and \ref{thm:main}.
\section{Preliminaries}\label{sec:parts}

In this section, we introduce some tools which plays a fundamental role in the proofs of the main theorems. We begin with an identity for sufficiently regular solutions of 
\eqref{eq:helmholtz-laplace}.
\begin{lemma}\label{lem:identity-laplace}
  Let $\Omega\subset\R^n$ be any Lipshitz domain, and let $\phi(|x|),\varphi(|x|):\R^n\to\R$ be two
  sufficiently regular radial multipliers. Then, any solution $u\in\mathcal H^1(\Omega)
  \cap\mathcal H^{\frac32+}_{\text{loc}}$ of 
  \eqref{eq:helmholtz-laplace}
  satisfies
  \begin{align}\label{eq:identity-laplace}
    & \int_\Omega\nabla_AuD^2\phi\overline{\nabla_Au}dx
    +\int_\Omega\varphi\left|\nabla_Au\right|^2dx
    -\frac14\int_\Omega\Delta\left(\Delta\phi
    +2\varphi\right)|u|^2dx
    \\
    & +\int_\Omega\left[\varphi V-\frac12\phi'(\partial_rV)\right]|u|^2dx
    -\Im\int_\Omega\phi'uB_\tau\cdot
    \overline{\nabla_Au}dx
    -k^2\int_\Omega\varphi|u|^2dx
    \nonumber
    \\
    &
    +\frac14\int_{\partial\Omega}
    |u|^2\left(\nabla\Delta\phi+2\nabla\varphi\right)
    \cdot\eta d\sigma(x)
    -\frac{k^2}{2}\int_{\partial\Omega}
    |u|^2(\nabla\phi\cdot\eta)d\sigma(x)
    \nonumber
    \\
    &
    +\frac12\int_{\partial\Omega}
    |u|^2V(\nabla\phi\cdot\eta)d\sigma(x)
    +\frac12\int_{\partial\Omega}
    |\nabla_A u|^2(\nabla\phi\cdot\eta)d\sigma(x)
    \nonumber
    \\
    &
    -\frac12\Re\int_{\partial\Omega}
    (\nabla_A u\cdot \eta)\overline u(\Delta\phi+2\varphi)d\sigma(x)
    -\Re\int_{\partial\Omega}
    (\nabla u\cdot\eta)(\nabla\phi\cdot\overline{\nabla_A u})d\sigma(x)
    \nonumber
    \\
    & = -\Re\int_\Omega
    f\left(\nabla\phi\cdot\overline{\nabla_Au}
    +\frac12(\Delta\phi)\overline{u}+\varphi\overline u\right)dx
    \mp\epsilon\Im\int_\Omega u\nabla\phi\cdot\overline{\nabla_A u}dx,
    \nonumber
  \end{align}
  where $\eta$ denotes the outgoing unit normal to the boundary $\partial\Omega$ and $d\sigma$ the surface measure on $\partial\Omega$;
   moreover, $D^2\phi,\Delta^2\phi$ denote, respectively, the Hessian and
  the bi-Laplacian of $\phi$, while $\partial_r V$ is the radial derivative of $V$ and $B_\tau$ is as in Definition
  \ref{def.bitau}.
\end{lemma}
The proof of the previous lemma is analogous to the one of Lemma 2.1 in \cite{F}, in which $\Omega=\R^n$ (see also \cite{pv} for the case $A\equiv0$, $\Omega=\R^n$).
Indeed, identity \eqref{eq:identity-laplace} turns out formally by summing up the two identities obtained as follows: 
\begin{enumerate}
  \item
  multiply equation \eqref{eq:helmholtz-laplace} by $\nabla\phi\cdot\overline{\nabla_Au}+
  \frac12(\Delta\phi)\overline u$, take the resulting real parts and integrate in $\Omega$;
  \item
  multiply equation by $\varphi\overline u$, take the resulting real parts and integrate in $\Omega$.
\end{enumerate}
We omit straightforward details, one should only follow carefully the boundary terms in $\partial\Omega$ which do not appear in \cite{F}.

See also remark 2.1 in \cite{F}, for a discussion on the regularity which is needed on $u,A,V$ in order to justify the above mentioned integration by parts.

\subsection{Choice of the multipliers}
\label{subsec:multipliers}
We now introduce some explicit multipliers; in this explicit form, they have been obtained by slightly modifying the ones used in 
\cite{F}, inspired to \cite{pv}.
In the following, we denote by $r=|x|$; in dimension $n\geq3$, let us define
 \begin{equation*}
  \phi_0(r)=\int_0^r\phi_0'(s)\,ds,
\end{equation*}
where
\begin{equation*}
  \phi'_0=\phi'_0(r)=
  \begin{cases}
    M+\frac{(n-1+\alpha)}{2n}r,
    \qquad
    r\leq1
    \\
    M+\frac12-\frac1{2nr^{n-1}}+\frac\alpha{2n},
    \qquad
    r>1,
  \end{cases}
\end{equation*}
and $M,\alpha>0$ are arbitrary constants. Observe that $\phi_0$ is a continuous function. By
scaling we define
\begin{equation*}
  \phi_R(r)=R\phi_0\left(\frac rR\right),
\end{equation*}
and by direct computations we obtain
\begin{equation}\label{eq:fi14d}
  \phi'_R=\phi'_0\left(\frac rR\right)=
  \begin{cases}
    M+\frac{(n-1+\alpha)}{2n}\cdot\frac rR,
    \qquad
    r\leq R
    \\
    M+\frac12-\frac{R^{n-1}}{2nr^{n-1}}+\frac{\alpha}
    {2n},
    \qquad
    r>R,
  \end{cases}
\end{equation}
\begin{equation}\label{eq:fi24d}
  \phi''_R=
  \begin{cases}
    \frac{n-1}{2nR}+\frac\alpha{2nR},
    \qquad
    r\leq R
    \\
    \frac{n-1}{2n}\cdot\frac{R^{n-1}}{r^n},
    \qquad
    r>R.
  \end{cases}
\end{equation}
Moreover we have
\begin{equation}\label{eq:lapphi}
  \Delta\phi_R(r)=
  \begin{cases}
    \frac{n-1}{2R}+\frac{M(n-1)}{r}+
    \frac{\alpha}{2R},
    \qquad
    r\leq R
    \\
    \frac{(2M+1)(n-1)}{2r}
    +\frac{(n-1)\alpha}{2nr},
    \qquad
    r>R,
  \end{cases}
\end{equation}
which defines a discontinuous function.
Now, let us define
\begin{equation}\label{eq:varfi}
  \varphi_R(r)=
  \begin{cases}
    -\frac\alpha{4nR}
    \quad
    r<R
    \\
    0
    \quad
    r\geq R,
  \end{cases}
\end{equation}
for any $R>0$.
We hence have
\begin{equation}\label{eq:deltafi4D}
  \Delta\phi_R(r)+2\varphi_R(r)=
  \begin{cases}
    \frac{n-1}{2R}+\frac{M(n-1)}{r}+\frac
    {(n-1)\alpha}{2nR},
    \qquad
    r\leq R
    \\
    \frac{(2M+1)(n-1)}{2r}
    +\frac{(n-1)\alpha}{2nr},
    \qquad\quad\
    r>R,
  \end{cases}
\end{equation}
which now defines a continuous function.
In view of identity \eqref{eq:identity-laplace}, we need to compute $\Delta\left(\Delta\phi_R+2\varphi_R\right)$. In dimension $n=3$ we have
\begin{equation}\label{eq:fibi}
  \Delta\left(\Delta\phi_R+2\varphi_R\right)
  =-4\pi M\delta_{x=0}
  -\frac{\alpha+3}{3R^2}\delta_{|x|=R},
\end{equation}
while in dimension $n\geq4$
the result is given by
\begin{align}\label{eq:bifi4D}
  \Delta\left(\Delta\phi_R+2\varphi_R\right)=
  &
  -\frac{(n-1)(n+\alpha)}{2nR^2}\delta_{|x|=R}
  -M\frac{(n-1)(n-3)}{r^3}\chi_{[0,R]}
  \\
  &
  -\left(M+\frac12\right)\frac{(n-1)(n-3)(2n+\alpha)}
  {2nr^3}
  \chi_{(R,+\infty)},
  \nonumber
\end{align}
in the distributional sense, where $\chi$ denotes the characteristic
function.

Finally, we denote by
\begin{equation}\label{eq:translate}
  \phi_{R,y}(x)=\phi_R(|x-y|),
  \qquad
  \varphi_{R,y}(x)=\varphi_R(|x-y|).
\end{equation}

\section{Proof of Theorems \ref{thm:lap3}, \ref{thm:main}}\label{sec:proof}
We are now ready to prove our main theorems. We start with the case of Dirichlet boundary conditions.

\subsection{Proof of Theorem \ref{thm:lap3}}\label{subsec:lap}
First notice that, by a density argument, it is not restrictive to assume that $f\in L^2(E)$.
Let us start by estimating the RHS of identity \eqref{eq:identity-laplace}, with the choice of $\phi,\varphi$ given in the previous section.
It is sufficient to prove the Theorem in the case $y=0$, choosing $\phi_R,\varphi_R$. The proof in the general case $y\neq0$ is completely analogous, modulo translating the multipliers as in \eqref{eq:translate}.
\subsubsection{Estimate of the RHS in
\eqref{eq:identity-laplace}}\label{subsubsec1}
Let us now put $\Omega=E$ in identity \eqref{eq:identity-laplace}; in the following, we extend to 0 $f$ and $u$ outside $E$ and consider integrals on the whole space.
By \eqref{eq:fi14d} and the Cauchy-Schwartz inequality, we have
\begin{align}\label{eq:estf}
  & \left|\int f\nabla\phi_R\cdot\overline{\nabla_A u}\right|
  \leq C\sum_{j\in\mathbb
  Z}\int_{C(j)}|f|\cdot|\nabla_A u|
  \\
  &\ \ \ \ \ \ \ \
  \leq C\sum_{j\in\mathbb Z}\left(2^{-j-1}\int_{C(j)}|\nabla_A u|^2\right)^{\frac12}
  \left(2^{j+1}\int_{C(j)}|f|^2\right)^{\frac12}
  \nonumber
  \\
  &\ \ \ \ \ \ \ \
  \leq C\left(\sup_{R>0}\frac1R\int_{|x|\leq R}|\nabla_A u|^2\right)^{\frac12}
  \sum_{j\in\mathbb Z}\left(2^{j+1}\int_{C(j)}|f|^2\right)^{\frac12}
  \nonumber
  \\
  &\ \ \ \ \ \ \ \
  \leq\gamma\sup_{R>0}\frac1R\int_{|x|\leq R}|\nabla_Au|^2
  +C(\gamma)N(f)^2,
  \nonumber
\end{align}
 for$\gamma,C(\gamma)>0$. Analogously, by \eqref{eq:deltafi4D}
\begin{align}\label{eq:estf2}
  & \left|\int f(\frac12\Delta\phi_R+\varphi_R)\overline{u}\right|
  \leq C\sum_{j\in\mathbb
  Z}\int_{C(j)}|f|\cdot\frac{|u|}{|x|}
  \\
  &\ \ \ \ \ \ \ \
  \leq C\sum_{j\in\mathbb Z}\left(2^{-j}\int_{C(j)}\frac{|u|^2}{|x|^2}\right)^{\frac12}
  \left(2^{j}\int_{C(j)}|f|^2\right)^{\frac12}
  \nonumber
  \\
  &\ \ \ \ \ \ \ \
  \leq C\left(\sup_{R>0}\frac1{R^2}\int_{|x|=R}|u|^2d\sigma\right)^{\frac12}
  \sum_{j\in\mathbb Z}\left(2^{j}\int_{C(j)}|f|^2\right)^{\frac12}
  \nonumber
  \\
  &\ \ \ \ \ \ \ \
  \leq\gamma\sup_{R>0}\frac1{R^2}\int_{|x|=R}|u|^2
  d\sigma
  +C(\gamma)N(f)^2.
  \nonumber
\end{align}
It remains now to estimate the last term at the RHS of
\eqref{eq:identity-laplace}. Multiplying \eqref{eq:diric-laplace}
by $u$ in $L^2$ and taking the resulting imaginary parts,
we obtain
\begin{equation}\label{eq:epsilon}
  |\epsilon|\int|u|^2dx\leq\int|fu|dx;
\end{equation}
analogously, taking the real parts we have
\begin{equation*}
  \int|\nabla_Au|^2=-\int V|u|^2+k^2\int|u|^2-\Re\int f\overline
  u.
\end{equation*}
Hence by assumption \eqref{eq:hardy}  we can estimate
\begin{align}\label{eq:estepsilon}
  & \left|\epsilon\int u\nabla\phi_R\cdot\overline{\nabla_Au}\right|
  \leq
  C|\epsilon|^{1/2}\left(\int(V_-+k^2)|u|^2+\int|fu|
  \right)^{\frac12}
  \left(\int|fu|\right)^{\frac12}
  \\
  &\ \ \ \ \ \ \ \
  \leq C|\epsilon|^{1/2}\int|fu|+
  C\left(|\epsilon|(k^2+\|V_-\|_\infty)\int|fu|\int
  |u|^2
  \right)^{\frac12}
  \nonumber
  \\
  &\ \ \ \ \ \ \ \
  \leq C(|\epsilon|+k^2+\|V_-\|_\infty)
  ^{\frac12}\int|fu|
  \nonumber
  \\
  &\ \ \ \ \ \ \ \
  \leq
  C(|\epsilon|+k^2+\|V_-\|_\infty)^{\frac12}
  \left(\sup_{R>0}\frac{k^2}R\int_{|x|\leq R}|u|^2\right)^{\frac12}
  \cdot N\left(\frac{f}{|k|}\right)
  \nonumber
  \\
  &\ \ \ \ \ \ \ \
  \leq\gamma\sup_{R>0}\frac{k^2}R\int_{|x|\leq R}|u|^2
  +C(\gamma)(|\epsilon|+k^2+\|V_-\|_\infty)
  \left(N\left(\frac{f}{|k|}\right)
  \right)^2,
  \nonumber
\end{align}
for $\gamma,C(\gamma)>0$. Recollecting \eqref{eq:estf},
\eqref{eq:estf2} and \eqref{eq:estepsilon}, for the right-hand side
of \eqref{eq:identity-laplace} we have
\begin{align}\label{eq:estRHS}
  & \left|-\Re\int_E
    f\left(\nabla\phi_R\cdot\overline{\nabla_Au}
    +\frac12(\Delta\phi_R)\overline{u}
    +\varphi_R\overline u\right)dx
    \mp\epsilon\Im\int_E u\nabla\phi_R\cdot\overline{\nabla_A u}dx\right|
  \\
  &\ \
  \leq\gamma\sup_{R>0}\left(\frac1R\int_{E\cap B(R)}
  |\nabla_A u|^2+\frac{k^2}R\int_{E\cap B(R)}|u|^2
  +\frac1{R^2}
  \int_{S(R)}|u|^2d\sigma\right)
  \nonumber
  \\
  &\ \ \ \ \
  +C(\gamma)\left[N(f)^2+(|\epsilon|+k^2)
  \left(N\left(\frac{f}{|k|}
  \right)\right)^2\right],
  \nonumber
\end{align}
for an arbitrary $\gamma>0$.

We now give another estimate of right-hand side, which will be fundamental in the sequel.
Recall again that
\begin{equation*}
  \nabla\phi_R\in L^\infty,
  \qquad
  \Delta\phi_R\sim\varphi_R\leq\frac{C}{|x|};
\end{equation*}
hence, using the magnetic Hardy inequality
\begin{equation*}
  \int\frac{|u|^2}{|x|^2}\,dx\leq\left(\frac{2}{n-2}\right)^2\int|\nabla_A u|^2\,dx
  \qquad
  (n\geq3)
  \quad
\end{equation*}
(see \cite{FV} for the proof), one obtains immediately that
\begin{align}\label{eq:finiti}
  & \left|-\Re\int_E
    f\left(\nabla\phi_R\cdot\overline{\nabla_Au}
    +\frac12(\Delta\phi_R)\overline{u}
    +\varphi_R\overline u\right)dx
    \mp\epsilon\Im\int_E u\nabla\phi_R\cdot\overline{\nabla_A u}dx\right|
  \\
  &\ \
  \leq C\left(\|f\|_{L^2}+\|\nabla_Au\|_{L^2}+\|u\|_{L^2}\right)<\infty,
  \nonumber
\end{align}
by standard elliptic theory, 
since we have assumed a priori that $f\in L^2(E)$.
 Here we identified $u$ and $f$ with their trivial extensions to the whole space.

Our next step is to prove the positivity of the left-hand side of
\eqref{eq:identity-laplace}.

\subsubsection{Positivity of the LHS of \eqref{eq:identity-laplace}}\label{sec:positivity}
Let us start by estimating the terms at the left-hand side of \eqref{eq:identity-laplace} involving solid integrals in $E=\Omega$. Again, by extending to 0 the solution $u$ outside $E$, we can consider integrals on the whole space.
The argument is slightly different in the 3D-case and in the higher dimensional case.

{\bf 3D-case.}
Let us start by the first term. Since $\phi_R$ is radial, the following formula holds:
\begin{equation}\label{eq:tang}
  \nabla_AuD^2\phi_R\overline{\nabla_Au}=
  \phi_R''|\nabla_A^ru|^2+\frac{\phi'_R}{|x|}
  |\nabla_A^\bot u|^2,
\end{equation}
where $\nabla_A^ru=\nabla_Au\cdot x/|x|$ and $|\nabla_A^\bot u|$ the
modulus of any tangent vector to the unit sphere in $x/|x|$, such that
\begin{equation*}
  \nabla_A^\bot u\cdot\nabla_A^ru=0,
  \qquad
  |\nabla_A^\bot u|^2=|\nabla_Au|^2-|\nabla_A^ru|^2.
\end{equation*}

By \eqref{eq:tang}, \eqref{eq:fi14d}, \eqref{eq:fi24d}, and \eqref{eq:fibi}
we can hence estimate
\begin{align}\label{eq:LHS1}
  &
  \int\nabla_AuD^2\phi_R\overline{\nabla_Au}
  +\int\varphi_R|\nabla_Au|^2
  -\frac14\int\Delta\left(\Delta\phi_R+2\varphi_R\right)|u|^2\,dx
  \\
  & \ \ \
  \geq
  M\int\frac{|\nabla_A^\bot u|^2}{|x|}
  +\frac CR\int_{|x|\leq
  R}|\nabla_Au|^2
  +\frac{C}{R^2}\int_{|x|=R}|u|^2d\sigma+C|u(0)|^2,
  \nonumber
\end{align}
for some $C>0$.

We now pass to the terms
containing $\partial_rV$ and $B_\tau$. By \eqref{eq:fi14d} we obtain
\begin{align}\label{eq:LHS3}
  & -\frac12\int\phi'_R(\partial_rV)|u|^2
  \geq
  C\int(\partial_rV)_-|u|^2-C
  \int(\partial_rV)_+|u|^2
  \\
  & \geq C\int(\partial_rV)_-|u|^2-C
  \int_0^\infty d\rho\int_{|x|=\rho}(\partial_rV)_+|u|^2d\sigma
  \nonumber
  \\
  & \geq C\int(\partial_rV)_-|u|^2-C
  \sup_{R>0}\left(\frac1{R^2}\int_{|x|=R}|u|^2
  d\sigma\right)
  \||\cdot|^2(\partial_rV)_+\|_{L^1_rL^\infty(S_r)},
  \nonumber
\end{align}
for some $C>0$.
Moreover we have
\begin{align}\label{eq:LHS33}
  & \int\varphi_RV|u|^2
  =
  \frac{\alpha}{4nR}\int_{|x|\leq R}
  V_-|u|^2
  -\frac{\alpha}{4nR}\int_{|x|\leq R}V_+|u|^2
  \\
  & \ \ \
  \geq\frac CR\int_{|x|\leq R} V_-|u|^2-C\int|x|^{-1}
  V_+|u|^2
  \nonumber
  \\
  &\ \ \
  \geq \frac CR\int_{|x|\leq R} V_-|u|^2
  -C\sup_{R>0}\left(\frac1{R^2}
  \int_{|x|=R}|u|^2
  d\sigma\right)\|\,|\cdot|
  V_+\|_{L^1_rL^\infty(S_r)},
  \nonumber
\end{align}
for some $C>0$.
For the term containing $B_\tau$ observe that
\begin{equation*}
  \left|B_\tau\cdot\nabla_Au\right|=
  |B_\tau||\nabla_A^\bot u|,
\end{equation*}
since $B_\tau$ is a tangential vector to the sphere; as a consequence we get
\begin{align}\label{eq:LHS4}
  & -\Im\int_{\R^n}u\phi'_RB_\tau\cdot
  \overline{\nabla_Au}\,dx
  \geq-C\int_{\R^n}|u|\cdot|B_\tau|\cdot
  |\nabla_A^\bot
  u|\,dx
  \\
  &\ \ \
  \geq-C\left(\int\frac{|\nabla_A^\bot
  u|^2}{|x|}\right)^{\frac12}
  \left(\int_0^{+\infty}d\rho
  \int_{|x|=\rho}|x|\cdot|u|^2\cdot|B_\tau|^2
  d\sigma\right)^{\frac12}
  \nonumber
  \\
  &\ \ \
  \geq-C\left(\int\frac{|\nabla_A^\bot
  u|^2}{|x|}\right)^{\frac12}
  \left(\sup_{R>0}\frac1{R^2}\int_{|x|=R}|u|^2
  d\sigma\right)^{\frac12}\||\cdot|^{\frac32}
  B_\tau\|_{L^2_rL^\infty(S_r)}.
  \nonumber
\end{align}
Let us introduce the following notations:
\begin{equation*}
  a:=\left(\int\frac{|\nabla_A^\bot
  u|^2}{|x|}\right)^{\frac12};
  \qquad
  b:=\left(\sup_{R>0}\frac1{R^2}\int_{|x|=R}|u|^2
  d\sigma\right)^{\frac12}.
\end{equation*}
Moreover, according to assumption \eqref{eq:smallness3}, we denote
\begin{equation*}
  C_1:=\||\cdot|^{\frac32}B_\tau\|_{L^2_rL^\infty(S_r)};
\end{equation*}
\begin{equation*}
  C_2:=\||\cdot|^2(\partial_rV)_+\|_{L^1_rL^\infty(S_r)};
\end{equation*}
\begin{equation*}
  C_3:=\|\,|\cdot|V_+\|_{L^1_rL^\infty(S_r)}.
\end{equation*}
We are ready now to sum \eqref{eq:LHS1},
\eqref{eq:LHS3}, \eqref{eq:LHS33}, and \eqref{eq:LHS4}. Since $R$ is arbitrary, we can take the supremum over $R$ in
\eqref{eq:LHS1}; it turns out that
\begin{align}\label{eq:LHS5}
  & \int\nabla_AuD^2\phi_R\overline{\nabla_Au}
  +\int\varphi_R|\nabla_Au|^2
  -\frac14\int\Delta\left(\Delta\phi_R+2\varphi_R\right)|u|^2\,dx
  \\
  &
  +\int\left(\varphi_RV
  -\frac12\phi'_R(\partial_rV)\right)|u|^2
  -\Im\int_{\R^n}u\phi'_RB_\tau\cdot
  \overline{\nabla_Au}
  \nonumber
  \\
  &
  \geq C\sup_{R>0}\frac1R\int_{|x|\leq
  R}|\nabla_Au|^2
  +C\sup_{R>0}\frac1R\int_{|x|\leq R}V_-|u|^2
  +C\int(\partial_rV)_-|u|^2
  \nonumber
  \\
  &\ \ \
  +Ma^2
  -CC_1ab
  +[C-CC_2-CC_3]b^2.
  \nonumber
  \\
  &
  \geq C\left(a^2b^2+\sup_{R>0}\frac1R\int_{|x|\leq
  R}|\nabla_Au|^2
  +\sup_{R>0}\frac1R\int_{|x|\leq R}V_-|u|^2
  +\int(\partial_rV)_-|u|^2\right),
  \nonumber
\end{align}
by taking $\delta$ small enough in condition \eqref{eq:smallness3}, for some constant $C>0$.
Since
\begin{equation}\label{eq:k}
  -k^2\int\varphi_R|u|^2
  \geq
  C\frac{k^2}{R}\int_{|x|\leq R}|u|^2,
\end{equation}
by \eqref{eq:LHS5} and \eqref{eq:k} we conclude that
\begin{align}\label{eq:LHS3d}
  & \int_E\nabla_AuD^2\phi_R\overline{\nabla_Au}dx
    +\int_E\varphi_R\left|\nabla_Au\right|^2dx
    -\frac14\int_E\Delta\left(\Delta\phi_R
    +2\varphi_R\right)|u|^2dx
    \\
    & +\int_E\left[\varphi_R V-\frac12\phi'_R(\partial_rV)\right]|u|^2dx
    -\Im\int_E\phi'_RuB_\tau\cdot
    \overline{\nabla_Au}dx
    -k^2\int_E\varphi_R|u|^2dx
    \nonumber
    \\
    &
    \geq
    C\sup_{R>0}\left(
    \frac1R\int_{E\cap B(R)}\left(|\nabla_Au|^2
    +|k|^2\cdot|u|^2\right)
    +\frac1{R^2}\int_{E\cap S(R)}|u|^2
    +\int\frac{|\nabla_A^\bot u|^2}{|x|}\right)
    \nonumber
    \\
    &
    \ \ \
    +C\sup_{R>0}\frac1R\int_{E\cap B(R)}V_-|u|^2
    +C\int_{E}(\partial_rV)_-|u|^2+C|u(0)|^2,
\end{align}
for some $C>0$.
In particular, \eqref{eq:identity-laplace}, \eqref{eq:finiti}, and \eqref{eq:LHS3d} we obtain that
\begin{equation}\label{eq:finitifiniti}
  \sup_{R>0}\left(
    \frac1R\int_{E\cap B(R)}\left(|\nabla_Au|^2
    +|k|^2\cdot|u|^2\right)
    +\frac1{R^2}\int_{E\cap S(R)}|u|^2
    \right)<\infty.
\end{equation}

We now pass to the boundary terms at the left-hand side of \eqref{eq:identity-laplace}.

Due to the
Dirichlet boundary conditions, all the terms containing $u$ disappear; moreover notice that
$u\vert_{\partial\sigma}\equiv0\Rightarrow\nabla^\tau u\vert_{\partial\sigma}\equiv0$, and consequently
\begin{equation*}
  \nabla_Au\vert_{\partial E} \equiv \nabla u\vert_{\partial E}\equiv \nabla^\eta u\vert_{\partial E}
  \equiv \nabla_A^\eta u\vert_{\partial E}.
\end{equation*}
Hence
\begin{equation}\label{eq:boundary1}
  \frac12\int_{\partial E}
    |\nabla_A u|^2(\nabla\phi_R\cdot\eta)d\sigma(x)
  =
  \frac12\int_{\partial E}
    |\nabla_A^\eta u|^2\phi_R'
    \frac{x}{|x|}\cdot\eta.
\end{equation}
Moreover, decomposing
\begin{equation*}
  \frac{x}{|x|}=\left(\frac{x}{|x|}\cdot\eta\right)
  \eta+\nu,
\end{equation*}
where $\nu\cdot\eta=0$, we obtain
\begin{equation}\label{eq:boundary2}
  -\Re\int_{\partial E}
    (\nabla u\cdot\eta)(\nabla\phi\cdot\overline{\nabla_A u})
    =
    -\int_{\partial E}|\nabla_A^\eta u|^2
    \left(\frac{x}{|x|}\cdot\eta\right)\phi_R'
\end{equation}
Summing up \eqref{eq:boundary1} and \eqref{eq:boundary2} we get
\begin{align}\label{eq:boundary}
  &
  \frac12\int_{\partial E}
    |\nabla_A u|^2(\nabla\phi_R\cdot\eta)d\sigma(x)
    -\Re\int_{\partial E}
    (\nabla u\cdot\eta)(\nabla\phi\cdot\overline{\nabla_A u})d\sigma(x)
    \\
  &
  =-\frac12\int_{\partial E}|\nabla_A^\eta u|^2
    \left(\frac{x}{|x|}\cdot\eta\right)\phi_R'
  \nonumber
  \\
  &
  \geq
  C_1\int_{\partial E^-_0}|\nabla_A^\eta u|^2\left(-\frac{x}{|x|}\cdot\eta\right)d\sigma
  -C_2\int_{\partial E^+_0}|\nabla_A^\eta u|^2\left(\frac{x}{|x|}\cdot\eta\right)
  d\sigma,
\end{align}
for some $C_1,C_2>0$,
since $\phi'_R
\in[M,M+(n+\alpha)/2n]$ by \eqref{eq:fi14d}.
Estimate \eqref{eq:estimate1} now follows by \eqref{eq:estRHS}, \eqref{eq:LHS3d} and \eqref{eq:boundary}, up to choose $\gamma$ sufficiently small in \eqref{eq:estRHS}, and the proof is complete. Indeed, one should absorb the first three terms at the right-hand side of 
\eqref{eq:estRHS}, which are finite by \eqref{eq:finitifiniti}, for a sufficiently small $\gamma>0$.

We can now pass to the higher dimensional case.

{\bf Higher dimensional case.} The proof is slightly different from the 3D case. Also in this case we prove \eqref{eq:estimate1} in the case $y=0$; indeed the general case $y\neq0$ is completely analogous, modulo the translation in \eqref{eq:translate}.

For the right-hand side of \eqref{eq:identity-laplace} we follow using estimate
\eqref{eq:estRHS}. For the positivity of the left-hand side, let us consider the solid terms in \eqref{eq:identity-laplace}.
Arguing as in the previous section,
by \eqref{eq:tang}, \eqref{eq:fi14d}, \eqref{eq:fi24d}, and \eqref{eq:bifi4D}
we can estimate
\begin{align}\label{eq:LHS14}
  &
  \int\nabla_AuD^2\phi_R\overline{\nabla_Au}
  +\int\varphi_R|\nabla_Au|^2
  -\frac14\int\Delta\left(\Delta\phi_R+2\varphi_R\right)|u|^2
  \\
  & \ \ \
  \geq
  M\int\frac{|\nabla_A^\bot u|^2}{|x|}
  +\frac CR\int_{|x|\leq
  R}|\nabla_Au|^2
  +\frac{C}{R^2}\int_{|x|=R}|u|^2d\sigma
  +C\int\frac{|u|^2}{|x|^3},
  \nonumber
\end{align}
for some $C>0$. For the terms containing $V$ and $B_\tau$ we now estimate
\begin{align}\label{eq:LHS34D}
  & -\frac12\int\phi'_R(\partial_rV)|u|^2
  \geq
  C_1\int(\partial_rV)_-|u|^2
  -C_2
  \|\,|\cdot|^3(\partial_rV)_+\|
  _{L^\infty}\int\frac{|u|^2}{|x|^3}dx,
  \\
  &
  \int\varphi_RV|u|^2\geq\frac{C_1}{R}\int_{|x|\leq R}V_-|u|^2
  -C_2\|\,|\cdot|^{-2}V_+\|
  _{L^\infty}\int\frac{|u|^2}{|x|^3}dx,
  \label{eq:LHS34D2}
  \\
  & -\Im\int u\phi'_RB_\tau\cdot
  \overline{\nabla_Au}\,dx
  \geq-C\left(\int\frac{|\nabla_A^\bot
  u|^2}{|x|}\right)^{\frac12}\left(
  \int\frac{|u|^2}{|x|^3}\right)^{\frac12}
  \|\,|\cdot|^2B_\tau\|_{L^\infty}.
  \label{eq:LHS34D3}
\end{align}
Now we can sum \eqref{eq:LHS14},
\eqref{eq:LHS34D}, \eqref{eq:LHS34D2} and \eqref{eq:LHS34D3}, taking
the supremum over $R$; we denote by
\begin{equation*}
  a:=\left(\int\frac{|\nabla_A^\bot u|^2}{|x|}
  \right)^{\frac12};
  \qquad
  b:=\left(\int\frac{|u|^2}{|x|^3}\right)^{\frac12},
\end{equation*}
and according to assumption \eqref{eq:smallness4}
\begin{equation*}
  \|\,|\cdot|^2B_\tau\|_{L^\infty}=C_1,
  \quad
  \|\,|\cdot|^3(\partial_rV)_+\|_{L^\infty}= C_2,
  \quad
  \|\,|\cdot|^2V_+\|_{L^\infty}=C_3
\end{equation*}
We obtain
\begin{align}\label{eq:LHS54D}
  & \int\nabla_AuD^2\phi_R\overline{\nabla_Au}
  +\int\varphi_R|\nabla_Au|^2
  -\frac14\int\left(
  \Delta(\Delta\phi_R+2\varphi_R)\right)|u|^2
  \\
  &
  +\int\left(\varphi_RV-
  \frac12\phi'_R(\partial_rV)+\right)|u|^2
  -\Im\int u\phi'_RB_\tau\cdot
  \overline{\nabla_Au}
  \nonumber
  \\
  &
  \geq C\sup_{R>0}\left(\frac1R\int_{|x|\leq R}
  |\nabla_Au|^2+\frac1{R^2}
  \int_{|x|=R}|u|^2\right)
  \nonumber
  \\
  &
  \ \ \ +C\int(\partial_rV)_-|u|^2
  +\sup_{R>0}\frac CR\int_{|x|\leq R}V_-|u|^2
  \nonumber
  \\
  &
  \ \ \ +Ca^2-CC_1ab+(C-CC_2-CC_3)b^2,
  \nonumber
\end{align}
for some $C>0$. Finally, by condition \eqref{eq:smallness4} we have
\begin{equation*}
  Ca^2-CC_1ab+(C-CC_2-CC_3)b^2
  \geq C(a^2+b^2),
\end{equation*}
for $\delta$ sufficiently small; as a consequence, \eqref{eq:LHS54D} implies
\begin{align}\label{eq:LHSfinal4D}
  & \int_E\nabla_AuD^2\phi_R\overline{\nabla_Au}
  +\int_E\varphi_R|\nabla_Au|^2
  -\frac14\int_E
  \Delta(\Delta\phi_R+2\varphi_R)|u|^2
  \\
  &
  +\int_E\left(\varphi_RV-
  \frac12\phi'_R(\partial_rV)+\right)|u|^2
  -\Im\int_Eu\phi'_RB_\tau\cdot
  \overline{\nabla_Au}
  \nonumber
  \\
  &
  \geq C\sup_{R>0}\left(\frac1R\int_{E\cap B(R)}
  |\nabla_Au|^2+\frac1{R^2}
  \int_{E\cap S(R)}|u|^2\right)
  +C\int_E\left(\frac{|u|^2}{|x|^3}+
  \frac{|\nabla_A^\bot u|^2}{|x|}\right)
  \nonumber
  \\
  &
  \ \ \ +C\int_E(\partial_rV)_-|u|^2
  +\sup_{R>0}\frac CR\int_{E\cap B(R)}V_-|u|^2.
  \nonumber
\end{align}
Notice the additional estimate for the term involving $\int|u|^2/|x|^3$, which does not hold in dimension $n=3$. From now on, the proof is completely analogous to the 3D-case.

\subsection{Proof of Theorem \ref{thm:main}}
The scheme of the proof is completely analogous to the previous one. In particular, the control of the terms in \eqref{eq:identity-laplace} which involve the integrals in $E$ work exactly as in the previous proof. Hence we just need now to control all the boundary terms in identity 
\eqref{eq:identity-laplace}, in the general case in which no boundary conditions are assumed.

 First, by \eqref{eq:fi14d} and differentiating in \eqref{eq:deltafi4D}  we get
 \begin{align}\label{eq:222}
 &
  \frac14\int_{\partial E}|u|^2
 \nabla\left(\Delta\phi_R+2\varphi_R\right)
 \cdot\eta
 -\frac{k^2}{2}\int_{\partial E}|u|^2
  \left(\nabla\phi_R
  \cdot\eta\right)
 \\
 &
 \geq
  C_1\int_{\partial E}|u|^2
  \left(-\frac{x}{|x|}\cdot\eta\right)
  +C_2k^2
 \int_{\partial E}|u|^2
  \left(-\frac{x}{|x|}\cdot\eta\right)
  \nonumber
  \\
  &
  \geq
  C\beta(1+k^2)\int_{\partial E}|u|^2,
  \nonumber
\end{align}
for some constant $C>0$.
 For the term containing $V$ we easily estimate
\begin{equation}\label{eq:444}
  \frac12\int_{\partial E}|u|^2V\left(\nabla\phi_R\cdot\eta\right)
  \geq
  -C\|V\|_{L^\infty(\partial E)}\int_{\partial E}|u|^2,
\end{equation}
for some $C>0$. 

Notice that, since $0\notin\partial\Omega$, then by \eqref{eq:fi14d} and \eqref{eq:deltafi4D}
we have
\begin{equation*}
  0\leq\Delta\phi_R+2\varphi_R\leq C_1
  \qquad
  C_2\leq\phi'_R\leq C_3,
  \qquad
  (x\in\partial E)
\end{equation*}
for some constants $C_1,C_2,C_3>0$. Moreover recall that by assumption 
$\beta\leq -x/|x|\cdot\eta\leq1$;
consequently, we can estimate
\begin{align*}
  &
  -\frac12\Re\int_{\partial\Omega}
    (\nabla_A u\cdot \eta)\overline u(\Delta\phi_R+2\varphi_R)
    \\
    &
    \geq -C\int_{\partial E}|\nabla_A u\cdot \eta|\cdot|u|\geq
    -\frac C\beta\sup_{x\in\partial E}\frac{1}{\phi'_R}\int_{\partial E}|\nabla_A u\cdot \eta|\cdot|u|
    \phi'_R\left(-\frac{x}{|x|}\cdot\eta\right),
\end{align*}
for some $C>0$. Hence, by Cauchy-Schwartz and the elementary inequality 
$ab\leq\epsilon a^2+\frac1\epsilon b^2$, for any $\epsilon\in(0,1)$, we obtain
\begin{align}\label{eq:666}
&
  -\frac12\Re\int_{\partial\Omega}
    (\nabla_A u\cdot \eta)\overline u(\Delta\phi_R+2\varphi_R)
 \\
    & \ \ \geq
    -C\epsilon\int_{\partial E}|\nabla_A^\eta u|^2\phi'_R\left(-\frac{x}{|x|}\cdot\eta\right)
    -\frac C{\epsilon}\int_{\partial E}|u|^2\phi'_R\left(-\frac{x}{|x|}\cdot\eta\right),
    \nonumber
\end{align}
for some $C>0$ and any $\epsilon>0$. For the last term, first write
\begin{align}\label{eq:first}
&
  -\Re\int_{\partial E}
    (\nabla u\cdot\eta)(\nabla\phi_R\cdot\overline{\nabla_A u})
    \\
    &
    =
    -\Re\int_{\partial E}
    (\nabla_Au\cdot\eta)(\nabla\phi_R\cdot\overline{\nabla_A u})
    +\Im\int_{\partial E}
    u(A\cdot\eta)(\nabla\phi_R\cdot\overline{\nabla_A u}).
    \nonumber
\end{align}
Hence, writing $\nabla\phi_R=\phi_R'\frac{x}{|x|}$ and decomposing
\begin{equation*}
  \frac{x}{|x|}=\left(\frac{x}{|x|}\cdot\eta\right)
  \eta+\nu,
\end{equation*}
where $\nu\cdot\eta=0$, we obtain
\begin{align*}
  &
  -\Re\int_{\partial E}
    (\nabla u\cdot\eta)(\nabla\phi_R\cdot\overline{\nabla_A u})
    \\
    & \ \ 
    =
    \int_{\partial E}|\nabla_A^\eta u|^2
    \left(-\frac{x}{|x|}\cdot\eta\right)\phi'_R
    -\Re\int_{\partial E}(\nabla_Au\cdot\eta)\left(\overline{\nabla_Au}\cdot\nu \right)\phi'_R
    \\
    & \ \ \ \ \ 
    +\Im\int_{\partial E}u(A\cdot\eta)\left(\overline{\nabla_Au}\cdot\eta\right)
    \left(\frac{x}{|x|}\cdot\eta\right)\phi'_R
    +
    \Im\int_{\partial E}u(A\cdot\eta)\left(\overline{\nabla_Au}\cdot\nu\right)
    \phi'_R;
\end{align*}
finally, arguing as above we conclude that
\begin{align}\label{eq:777}
  &
  -\Re\int_{\partial E}
    (\nabla u\cdot\eta)(\nabla\phi_R\cdot\overline{\nabla_A u})
    \\
    & \ \ 
    \geq
    \left[1-\epsilon\left(\frac1\beta+1\right)\right]
    \int_{\partial E}|\nabla_A^\eta u|^2
    \left(-\frac{x}{|x|}\cdot\eta\right)\phi'_R
    -\frac C\epsilon
    \int_{\partial E}|\nabla_A^\tau u|^2
    \nonumber
    \\
    & \ \ \ \ \ 
    -\frac C\epsilon\|A\|_{L^\infty(\partial E)}^2
    \int_{\partial E}|u|^2,
    \nonumber
 \end{align}
for some $C>0$ and any $\epsilon>0$. 
Estimate \eqref{eq:estimate1main} now follows by \eqref{eq:identity-laplace}, \eqref{eq:222}, \eqref{eq:444}, \eqref{eq:666}, \eqref{eq:777}, choosing $\epsilon>0$ small enough, depending on $\beta$.
Finally, if $k$ is large enough, the term in \eqref{eq:222}, which is positive, controls all the other terms containing $|u|^2$, and also \eqref{eq:estimate1main2} is proved.

%%%%%%%%%%BIBLIOGRAPHY%%%%%%%%%%%%%

\end{document}